\documentclass{article}
\usepackage{amsmath,amsfonts, amsthm, amscd, amssymb}
\usepackage[dvips]{epsfig,graphics}

\numberwithin{equation}{section}

\newtheorem{theorem}{Theorem}[section]
\newtheorem{proposition}[theorem]{Proposition}
\newtheorem{lemma}[theorem]{Lemma}
\newtheorem{corollary}[theorem]{Corollary}

\theoremstyle{definition}
\newtheorem{definition}[theorem]{Definition}
\newtheorem{example}[theorem]{Example}
\newtheorem{remark}[theorem]{Remark}
\newtheorem{question}[theorem]{Question}

\newcommand{\Z}{\mathbb{Z}}
\newcommand{\GR}{\rm{GR}}
\newcommand\D{{\Delta}}
\newcommand\s{{\sigma}}

\def\={\ {\buildrel \cdot \over =}\ }

\newcommand\bibname{\def\and{{\normalfont\rmfamily and }}%
                    \normalfont\scshape}%

\begin{document}

\title{Lehmer's Question, Knots and Surface Dynamics}
\author{Daniel S. Silver \& Susan G. Williams  \thanks{Both authors partially supported by NSF grant
DMS-0304971.}}
\maketitle

\begin{abstract}
Lehmer's question is equivalent to one about generalized growth rates of Lefschetz numbers of iterated pseudo-Anosov surface homeomorphisms. One need consider only homeomorphisms that arise as monodromies of fibered knots in lens spaces $L(n,1), n> 0$. Lehmer's question for Perron polynomials is equivalent to one about generalized growth rates of words under injective free group endomorphisms.  
\end{abstract}



\section{Introduction}
\label{intro} 
The {\it Mahler measure} $M(f)$ of a nonzero polynomial $f(x)= c_nx^n + \cdots + c_1 x + c_0\ (c_n \ne 0) \in {\mathbb C}[x]$  is  
$$M(f)= |c_n|\prod_{j=1}^n {\rm max}(|r_j|,1),$$
where $r_1, \ldots, r_n$ are the roots of $f$. Lehmer's question asks whether the Mahler measure of a polynomial with integer coefficients can be arbitrarily close but not equal to $1$. D.H. Lehmer posed the question in \cite{Leh33}. He could do no better than 1.17628..., a value that he achieved with the remarkable polynomial $L(x)= x^{10}+x^9-x^7-x^6-x^5-x^4-x^3+x+1$. No value closer but not equal to $1$ has yet been found despite extensive computer searchers \cite{Boy80,Boy89,Mo98,Ra94}. History and status of various searches can be viewed at M. Mossinghoff's web page www.cecm.sfu.ca/$\sim$mjm/Lehmer. 

Mahler measure can be defined, albeit in a different manner, for polynomials of more than one variable.
Lehmer's question turns out to be equivalent to the identical question for polynomials in any fixed number of variables (see \cite{Boy81}).

Lehmer's polynomial evaluated at $x=-t$ is the Alexander polynomial of the $(-2,3,7)$-pretzel knot, a fibered hyperbolic knot with striking properties  \cite{Ki97}. Motivated by this, the authors applied techniques of algebraic dynamics  \cite{Sc95} to provide a topological interpretation of the Mahler measure of Alexander polynomials of knots and links in terms of homology of finite abelian branched covering spaces  \cite{SW02}.  Many hyperbolic knots and links with simple complements have Alexander polynomials with remarkably small Mahler measure  \cite{SW04}. 

All of the knots and links that arise in  \cite{SW04} are fibered, which means that their complements are mapping tori of surface homeomorphisms. It is natural to ask whether Lehmer's question is equivalent to a question about surface homeomorphisms. We answer this affirmatively in section \ref{sec3}. In section \ref{sec4} we go further and show that Lehmer's question is equivalent to a question about generalized growth rates of Lefschetz numbers of iterated pseudo-Anosov homeomorphisms. The growth rates that we use are a generalization of the usual one, introduced by Kronecker and Hadamard. 

Section \ref{sec5} explores  free group endomorphisms, growth rates of words and connections with Lehmer's question. The final section examines connections with braids. 

We are grateful to Vasiliy Prokhorov for informing us about Hada\-mard's work. Mike Boyle, David Fried, Eriko Hironaka, Michael Kelly and Alexander Stoimenow also contributed helpful comments.  We thank the referee for suggestions and corrections.


\section{Growth rates of sequences} 
\label{sec2}
Let ${\bf a} = \{a_n\}_{n=1}^\infty$ be a sequence of complex numbers. For each positive integer $k$, we define the {\it $k$th growth rate} $\GR^{(k)}({\bf a})$ to be 
\[\limsup_{n\to \infty} \ \Biggr  \vert \det\ 
\left(
\begin{array}{c@{\enskip}c@{\enskip}c@{\enskip}c}
a_n & a_{n+1} & \cdots & a_{n+k-1} \\ a_{n+1} & a_{n+2} & \cdots & a_{n+k} \\ 
\vdots & & & \vdots \\ a_{n+k-1} & a_{n+k} & \cdots & a_{n+2k -2}
\end{array}\right) \Biggr\vert ^{1/n}.\]

We define the {\it $0$th growth rate} $\GR^{(0)}({\bf a})$ to be $1$. Notice that $\GR^{(1)}({\bf a})$ is the usual growth rate 
$$\limsup_{n\to \infty} \vert a_n \vert^{1/n},$$
which appears for example in  \cite{Bow78}. The general growth rates $\GR^{(k)}({\bf a})$ were studied by Kronecker and Hadamard in the late 19th century. The determinants $H_{n,k}$ above are known as {\it Hankel determinants}. The following proposition is a consequence of  \cite{Ha92} (see also p. 335 of \cite{Di57}).

\begin{proposition} The sequence ${\bf a}$ is the sequence of coefficients of a rational power series
\begin{equation}\label{eq2.1}
\sum_{n=0}^\infty a_n t^n = R(t)/\prod_{i=1}^d (1-\lambda_i t) 
\end{equation} 
if and only if for some $N$ the Hankel determinants $H_{n,k}$ vanish for all $k>d$ and $n\ge N$. 
In this case, $\GR^{(k)}({\bf a})=0$ if $k>d$.
If  no $1-\lambda_i t$ divides $R(t)$ and $|\lambda_1| \ge \cdots \ge |\lambda_d|$, then $\GR^{(1)}({\bf a})= |\lambda_1|,\  \GR^{(2)}({\bf a})= |\lambda_1 \lambda_2|, \cdots,\  
\GR^{(d)}({\bf a})= |\lambda_1 \lambda_2\cdots \lambda_d|.$
Consequently,  
$\max_{k}\GR^{(k)}({\bf a}) = M(f),$
where $f(t)=\prod_{i=1}^d (t- \lambda_i)$.
\label{hankel}\end{proposition}

By a {\it linear recurrence} we will mean a homogenous linear recurrence relation with constant coefficients, that is, a relation
$$a_{n+d} + c_{d-1} a_{n+d-1} + \cdots + c_1 a_{n+1} + c_0 a_n =0,$$
$n \ge 0$, with {\it characteristic polynomial} 
$$f(t) = t^d + c_{d-1} t^{d-1} + \cdots + c_1 t + c_0 = \prod_{i=1}^d (t-\lambda_i).$$   (An excellent reference on this topic is   \cite{EPSW03}.) The sequence ${\bf a} = \{a_n\}$ satisfies this linear recurrence if and only if it is the sequence of coefficients of the rational power series (\ref{eq2.1}), where $R(t) = b_{d-1} t^{d-1} + \cdots + b_1 t + b_0$ with
$b_k = \sum_{i=0}^k a_i c_{k-i}$.  The {\it minimal polynomial} of the sequence ${\bf a}$ is the unique minimal degree characteristic polynomial of a linear recurrence for ${\bf a}$.  Note that if the rational function  (\ref{eq2.1}) is in reduced form then the minimal polynomial of ${\bf a}$ is  $t^m\prod_{i=1}^d (t- \lambda_i)$ for some non-negative integer $m$.  Thus we have:

\begin{proposition}\label{prop2.2}If ${\bf a}=\{a_n\}_{n=1}^\infty$ is a linearly recurrent sequence with minimal polynomial $f(t)$ then 
$$\max_k \GR^{(k)}({\bf a}) = M(f).$$  \end{proposition}

The following result will be needed in Section \ref{sec4}.

\begin{proposition}\label{prop2.3}Suppose ${\bf a}=\{a_n\}_{n=1}^\infty$ and  ${\bf b}=\{b_n\}_{n=1}^\infty$ are linearly recurrent sequences with $|a_n-b_n|$ bounded.  Then the minimal polynomials of ${\bf a}$ and ${\bf b}$ have the same roots outside the unit circle, with the same multiplicities.  In particular, $\max_k \GR^{(k)}({\bf a}) =\max_k \GR^{(k)}({\bf b})$.
 \end{proposition}

\begin{proof} Set ${\bf c}=\{c_n\}=\{a_n-b_n\}$.  Then ${\bf c}$ is also linearly recurrent.  The power series $\sum_{n=0}^\infty c_n t^n$ converges on $|t|<1$.  Writing this power series as a rational function as in  (\ref{eq2.1}), we see that the roots of the minimal polynomial of  ${\bf c}$ must lie in $|t|\le1$.  From the general theory of linear recurrences, we can write $a_n=\sum p_i(n)\alpha_i^n$,  $b_n=\sum q_i(n)\beta_i^n$, $c_n=\sum r_i(n)\gamma_i^n$, where $\alpha_i$, $\beta_i$ and $\gamma_i$ range over the distinct roots of the minimal polynomials of ${\bf a}$, ${\bf b}$ and ${\bf c}$ respectively, and $p_i$, $q_i$ and $r_i$ are polynomials of degree one less than the multiplicity of the corresponding root.  Any set of functions on ${\mathbb N}$ of the form $n^k\lambda^n$ with distinct $(k,\lambda)\in {\mathbb N}\times{\mathbb C}$ is linearly independent.  Hence the terms in the sums $\sum p_i(n)\alpha_i^n$ and  $\sum q_i(n)\beta_i^n$ corresponding to roots of modulus greater than one must be identical.  The desired result is immediate. 
\end{proof}

\begin{remark}\label{rmk2.4} If we assume in addition to the hypothesis of Proposition \ref{prop2.3} that {\bf a} and {\bf b} are integer sequences, then the difference $a_n -b_n$ is eventually periodic (see Part VIII, Problem 158 of  \cite{PS76}). \end{remark}


\section{Fibered links}\label{sec3} Let $\ell = \ell_1 \cup \cdots \cup \ell_d \subset M^3$ be an oriented link in an oriented closed $3$-manifold with exterior $X = M\setminus {\rm int} N(\ell)$. Here $N(\ell)= \ell \times {\mathbb D}^2$
denotes a regular neighborhood of $\ell$. The link is {\it fibered} if the projection 
$\partial N(\ell) = \ell \times{\mathbb S}^1 \to {\mathbb S}^1$ extends to a locally trivial fibration of $X$. In this case, $X$ is homeomorphic to a {\it mapping torus} $S\times [0,1]/ \{(s,0)\sim (h(s),1)\}$, for 
some compact orientable surface ({\it fiber}) $S$  and homeomorphism (called a {\it monodromy}) $h: S \to S$, which is unique up to isotopy.

A monodromy $h$ induces an automorphism $h_*$ of $H_1(S; {\mathbb R})\cong {\mathbb R}^d$. The  characteristic polynomial ${\rm char}(h_*)= \det (h_* - tI)$ is a monic reciprocal integral polynomial of even degree that is well defined up to multiplication by $\pm t^n$. When $M$ is the $3$-sphere, it can be obtained from the classical multivariable Alexander polynomial $\D_\ell(t_1, \ldots, t_d)$ of $\ell$ by replacing each variable $t_i$ by $t$, and dividing by $(t-1)^{d-1}$ if $d>1$. Details can be found for example in  \cite{Kaw96}, Proposition 7.3.10. 

A surface automorphism $h:S\to S$ is {\it periodic} if $h^r$ is isotopic to the identity for some $r>0$. It is {\it reducible} if it is isotopic to an automorphism that leaves invariant some essential $1$-manifold of $S$. A theorem of W. Thurston  \cite{Th83} asserts that if $h$ is neither periodic nor reducible and $S$ is hyperbolic, then after a suitable isotopy, $h$ leaves invariant a pair of transverse singular measured foliations $({\cal F}^u, \mu^u), ( {\cal F}^s, \mu^s)$ such that
$$f({\cal F}^u, \mu^u)= ({\cal F}^u, \lambda\mu^u),\quad f({\cal F}^s, \mu^s)=({\cal F}^s, \lambda^{-1}\mu^s),$$
for some $\lambda>1$. 
Consequently, $h$ expands leaves of ${\cal F}^u$ by a factor of $\lambda$, and it contracts leaves of ${\cal F}^s$
by $\lambda^{-1}$. 
Such an automorphism is said to be {\it pseudo-Anosov}. Details for closed surfaces can be read in  \cite{BC88}. Modifications needed for surfaces with boundary can be found in  \cite{JG93} and  \cite{La79$'$}. Another result of Thurston  \cite{Th82} states that the mapping torus of a surface automoprhism $h:S\to S$ is a hyperbolic $3$-manifold if and only if $h$ is isotopic to a pseudo-Anosov homeomorphism. 

The following lemma is proved in section 5 of  \cite{BC88} for the case of a closed surface. The argument applies equally well to surfaces with connected boundary. For the reader's convenience, we repeat the argument with necessary modifications.

\begin{lemma}\label{lemma3.1} Assume that $h:S\to S$ is a homeomorphism of an compact connected orientable surface, closed or with connected boundary. If ${\rm char}(h_*)$ is not composite, cyclotomic or a polynomial in $t^r$ for any $r>1$, then $h$ is neither periodic nor reducible. 
\end{lemma}

\begin{proof} If $h$ is periodic, then $h_*^n =1$ for some $n>0$, which implies that ${\rm char}(h_*)$ is cyclotomic, a contradiction. 

Assume that $h$ is reducible. Then after isotopy, $h(C)=C$ for some essential $1$-submanifold $C$. The components of $C$ are disjoint simple closed curves, none homotopic to a point or $\partial S$, and no two homotopic to each other. We distinguish two cases: 

(1) Some component $C_1$ is non-separating. Orient $C_1$. It follows that $0\ne [C_1]\in H_1(S)$. For some $n>0$, $h_*^n[C_1]= [C_1]$. 
Hence some eigenvalue of $h_*$  is a root of unity, and  
${\rm char}(h_*)$ is cyclotomic, a contradiction. 

(2) Each component of $C$ separates $S$. There is a component $S_0$ of $S\setminus C$ with frontier equal to a single component of $C$. If $S$ has boundary, then we can assume that $S_0$ does not contain $\partial S$. Note that $H_1(S_0)$ must be nontrivial. Consider the least $r>0$ such that $h^r(S_0)=S_0$. Then 
$$H_1(S) = H_1(S_0) \oplus H_1(h(S_0)) \oplus \cdots \oplus H_1(h^{r-1}(S_0)) \oplus H_1(\hat G),$$
where $\hat G$, possibly empty, is obtained from $G=S\setminus (S_0 \cup h(S_0) \cup \cdots \cup h^{r-1}(S_0))$ by capping off boundary components other than $\partial S$. 

Since the summands $H_1(h^i(S_0))$ are nontrivial, irreducibility of 
${\rm char}(h_*)$ implies that $H_1(\hat G)$ is trivial. It follows that $r >1$, for otherwise $C$ consists of a single curve that is homotopic through $G$ to the boundary of $S$ or to a point if $\partial S$ is emtpy. Since the automorphism $h_*$ permutes the summands $H_1(h^i(S_0))$ cyclically,  ${\rm char}(h_*)$ has the form ${\rm det}(B-t^rI)$, for some matrix $B$. 
Hence ${\rm char}(h_*)$ is a polynomial in $t^r$, a contradiction. \end{proof}

The  Alexander polynomial $\D(k)(t)$ of a fibered knot $k\subset {\mathbb S}^3$ is an integral monic polynomial that satisfies the conditions: (1) $\D(k)(t)$ is  reciprocal; (2) $\D(k)(t)$ has even degree; (3) $\D(k)(1)=\pm 1$. Conversely, any such polynomial arises from a fibered knot. When investigating Lehmer's question, it suffices to consider integral monic polynomials satisfying conditions (1) and (2) (see  \cite{Sm71}). Unfortunately, it is not known if (3) can also be assumed.  In order to circumvent this problem, E. Hironaka suggested that one consider reduced Alexander polynomials of fibered links. 

By   \cite{Kan81},  for any integral monic polynomial $f(t)$ satisfying conditions (1) and (2), there exists a fibered $2$-component  link $\ell$ with $\D(\ell)(t, t)
\=\ f(t)$. (Here and throughout, $\=$ indicates equality up to multiplication by monomials $\pm t^n$.) The polynomial $\D(\ell)(t,t)$ is sometimes called the 
{\it reduced Alexander polynomial} of $\ell$, denoted by $\D_{\rm red}(\ell)(t)$. 

\begin{example}\label{ex3.2} Recall that $L(t)$ denotes Lehmer's polynomial
(see Introduction). It is well known that $L(-t)$ is the Alexander polynomial of the $(-2,3,7)$-pretzel knot, a fibered hyperbolic knot. Addition of a meridianal circle produces a  fibered $2$-component link $\ell$ with $\D_{\rm red}(\ell)(t)=L(-t)$. Kanenobu's general construction produces a much more complicated link $\ell$ with $\D_{\rm red}(\ell)(t) = L(-t)$. The general form of Kanenobu's links appears in Figure 1. The link is displayed as the boundary of a fiber surface $S$. A knot with Alexander polynomial $L(-t)$ results from the choice $n=5,\ p_1=-11,\ p_2=43,\ p_3=-69,\ p_4=36$ and $p_5=1$. \end{example}

\begin{figure}
\begin{center}
\bigskip
{\includegraphics[width=5in]{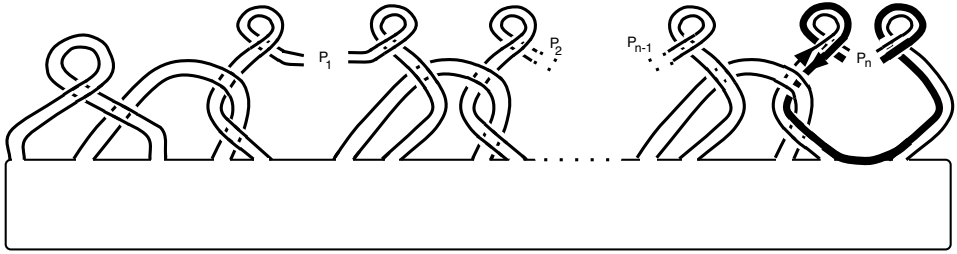}} 
\caption{Fibered link from Kanenobu's construction}
\label{kanenobu}
\end{center}
\end{figure}

In general, the links that Kanenobu constructs have the feature that at least one component is unknotted (shown in Figure \ref{kanenobu}). 

The definition of Alexander polynomial for knots in ${\mathbb S}^3$ extends for knots $k$ in rational homology $3$-spheres $M^3$  \cite{Tu01}. As in the classical case, when $k$ is fibered (that is, $M\setminus k$ fibers over the circle), the Alexander polynomial $\D(k)(t)$ coincides with ${\rm char}(h_*)$, where $h$ is a monodromy. It is well defined up to multiplication by $\pm t^i$.

Let $p, q$ be coprime integers. We recall that the {\it lens space} $L(p,q)$ is the union of two solid tori $V_1, V_2$ identified along their respective boundary tori $T_1, T_2$ by a homeomorphism from $V_2$ to $V_1$ that sends the homology  class of $\mu_2$ to that of $p\lambda_1 + q \mu_1$. 
Here $(\mu_i, \lambda_i)$ are respective meridian-longitude bases of $T_i,\ i=1,2$. It is a standard fact that $\pi_1(L(p,q))\cong \Z/p$. Also, $L(1,1)$ is the $3$-sphere ${\mathbb S}^3$. 
Additional background can be found, for example, in  \cite{Ro76}.

\begin{theorem}\label{thm3.3} Lehmer's question is equivalent to the following.
For any $\epsilon >0$, does there exist a fibered hyperbolic knot $k$ in a lens space $L(n,1)$, $n > 0,$ such that $1<M(\D(k)(t))<1+\epsilon$?
\end{theorem}

\begin{proof} We need consider only polynomials $f(t)$ which are not composite, cyclotomic or a polynomial in $t^r$ for any $r>1$. Kanenobu's construction produces a fibered $2$-component link $\ell= \ell_1\cup \ell_2$ with monodromy $h:S \to S$ such that ${\rm char}(h_*) = (t-1) f(t)$.
Capping off the boundary component of $S$  corresponding to an unknotted component, say $\ell_2$, produces a homeomorphism $\hat h: \hat S\to \hat S$ of a surface $\hat S$ with a single boundary component, and ${\rm char}(h_*)= f(t)$. By Lemma \ref{lemma3.1}, $\hat h$ is isotopic to a pseudo-Anosov homeomorphism.  

We prove that the mapping torus of $\hat h$ is the exterior of a knot $k$ in a lens space $L(n,1)$, where $n=f(1)$. (By  \cite{Ho58}, $n$ is the linking number of the two components of $\ell$.) 
Let $C_1, C_2$ denote the boundary components of $S$. Each is a longitude of the corresponding component $\ell_1, \ell_2$. Since $\ell_2$ is trivial, we may assume after isotopy that it is the standard unknot; let $V = \ell_2\times {\mathbb D}^2$ be a regular neighborhood of $\ell_2$. Let $V_1$ be the complementary solid torus in ${\mathbb S}^3$. Take $\lambda_1$ to be a meridian of $\ell_2$, $\mu_1$ a nontrivial curve in $\partial V_1$ that bounds a disk in $V$. 
   
The curve $\ell_1$ has winding number $n$ in $V_1$. Note that since $f(t)$ is irreducible and non-cyclotomic, $n$ is nonzero. The mapping torus of $\hat h$ is obtained from $V_1\setminus {\rm int} N(\ell_1)$ by attaching another solid torus $V_2= {\mathbb D}^2 \times {\mathbb S}^1$ with
$\partial {\mathbb D}^2 \times \theta = C_2 \times \theta$, for each $\theta \in {\Bbb S}^1$. It suffices to show that $C_2$ meets $\partial V_1$ in a curve homotopic to $\mu_1 + n \lambda_1$.  

It is clear that $C_2$ wraps once in the $\mu_1$ direction. The longitude $C_1$ of $\ell_1$ wraps $n$ times in the $\lambda_1$ direction. Since the fiber $S$ is a homology between $C_1$ and $C_2$, it follows immediately that $C_2$ also wraps $n$ times. This completes the argument that $\ell_1$ is a fibered knot in $L(n,1)$.

Since $L(-n, 1)$ is homeomorphic to $L(n,1)$, we may assume without loss of generality that $n>0$. 

Finally, we observe that since the monodromy $\hat h$ is pseudo-Anosov, the complement of $\ell_1$ is hyperbolic. \end{proof}
 
\begin{remark}\label{rem3.5}
(1) The knot $k$ can be obtained directly from $\ell$ by surgery on the unknotted component $\ell_2$.

(2) We can cap off the remaining boundary component of $\hat S$ and extend $\hat h$ to a homeomorphism $\bar h$ of a closed surface $\bar S$. The characteristic polynomials of $\hat h_*$ and $\bar h_*$ are equal.  Lemma 3.1 imples that Lehmer's question is equivalent to a question about pseudo-Anosov homeomorphisms of closed orientable surfaces.
\end{remark}


\section{Dynamics}\label{sec4} The periodic orbit structure of a surface homeomorphism contains much information about the dynamics of the  map. The fixed points of a pseudo-Anosov homeomorphism $h$ are isolated and hence can be counted. Since any iterate of $h$ is also pseudo-Anosov, the cardinality of the fixed point set ${\rm Fix}(h^n)$ is well defined for any integer $n$, and we denote it by 
$F_n$. One can regard ${\rm Fix}(h^n)$ as the  set of points of period $n$ under $h$. 

The {\it Lefschetz number} $L(h)$ of a pseudo-Anosov homeomorphism $h$ is the sum of the indices of its fixed points $p\in {\rm Fix}(h)$.  
(Lefschetz numbers are defined for more general maps.) If $p$ is an isolated fixed point in the interior of $S$, then its {\it index} is the degree of ${\rm id}-h: (S, S\setminus \{p\}) \to (S, S\setminus \{p\})$ . When the fixed point is contained on the boundary of $S$, then one computes the index by first embedding $S$ in Euclidean space as a neighborhood retract and then extending $h$. Details can be found in  \cite{Do80} (see also  \cite{Fe00}). 

By  \cite{FL79},
$2 \chi(S) -1 \le {\rm index}(p) \le 1.$
Moreover, the index at any non-prong singularity is $\pm 1$, and the number of all fixed points with index less than $-1$ is bounded below by $4 \chi(S)$  \cite{Ke97}.

We abbreviate the Lefschetz number $L(h^n)$ by $L_n$. In what follows, we consider the sequence ${\bf L} =\{L_n\}$. 

Lefschetz numbers  are homotopy invariants of $h$, and the well-known Lefschetz fixed point theorem implies that $L_n$ is equal to the alternating sum of homology traces:
\begin{eqnarray*}
L_n &= \displaystyle\sum_{i=0}^2 (-1)^i tr[h_*^n: H_i(S; {\mathbb R}) \to  H_i(S; {\mathbb R})]\\
  & = \begin{cases} 1-tr A^n, & \textrm{if}\ \partial S \ne \emptyset; \\
2- tr A^n, & \textrm{if}\ \partial S =\emptyset, \end{cases}\end{eqnarray*}
where $A$ is a matrix that represents $h_*:H_1(S, {\mathbb R})\to H_1(S, {\mathbb R})$. 

 Let $f(t) = \det (h_* - t I) $. Then $tr A^n = \lambda_1^n +\cdots + \lambda_d^n$, where $\lambda_1, \ldots,, \lambda_d$ are the roots (with possible repetition) of $f(t)$.  The sequence ${\bf L} = \{ L_n\}_{n=1}^\infty$ is linearly recurrent with minimal polynomial $p(t)=(t-1)\prod(t-\lambda_i)$, where the product ranges over distinct non-zero $\lambda_i$.   By Proposition \ref{prop2.2} we have:

\begin{proposition}\label{prop4.1} If $h$ is a pseudo-Anosov homeomorphism and $f(t) = \det (h_* - t I) $, then 
$$\max_k \GR^{(k)}({\bf L}) \le M(f),$$
with equality if $f$ has no repeated roots.\end{proposition}

Combining this with the proof of Theorem \ref{thm3.3} gives:
\begin{corollary}\label{cor4.2} Lehmer's question is equivalent to the following. For any $\epsilon >0$, does there exist a pseudo-Anosov homeomorphism $h: S\to S$ such that $\max_k \GR^{(k)}({\bf L})$ is contained in the interval $(1, 1+\epsilon)$?
\end{corollary}
One can assume that $h$ is a monodromy for a fibered 2-component link in ${\mathbb S}^3$ or a hyperbolic fibered knot in $L(n,1),\ n > 0$.  Alternatively, Remark \ref{rem3.5} implies that we can restrict our attention to homeomorphisms of closed orientable surfaces.

There exist pseudo-Anosov homeomorphisms $h$ such that $h_*$ is the identity map  \cite{Th88}, and consequently all growth rates $\GR^{(k)}{({\bf L})}$ are trivial. For such homeomorphisms, nontrivial growth rates are achieved by replacing Lefschetz numbers by the number $F_n$ of periodic points.
It is a consequence of  \cite{Man71} and  \cite{FS79} that ${\bf F}=\{F_n\}_{n=1}^\infty$ is linearly recurrent. We sketch the proof for the convenience of the reader. 

\begin{lemma}\label{lemma4.3} If $h: S \to S$ is a pseudo-Anosov homeomorphism, then ${\bf F}=\{F_n\}$ satisfies a linear recurrence. \end{lemma}

\begin{proof} By  \cite{FS79}, the homeomorphism $h:S\to S$ admits a Markov partition. Consequently, there exists an $N\times N$ matrix 
\begin{equation}\label{eq4.1} A = (a_{i,j})_{(i,j) \in N\times N} \end{equation}
with entries equal to $0$ or $1$ and a surjective map $p: \Omega \to S$, where 
$$\Omega = \{(x_n)_{n \in {\mathbb Z}} \mid a_{x_n, x_{n+1}} = 1,\ n \in {\mathbb Z}\}$$
such that $p \circ \sigma = h \circ p$, where $\sigma$ is the shift (to the left) of the sequence $(x_n)$ of symbols. The homeomorphism $\sigma$ is the {\it shift of finite type} with adjacency matrix $A$. For any $n$, the number of fixed points of $\sigma^n$ is equal to the trace of $A^n$ (see  \cite{BL68}, for example). However, $p$ need not be injective, and the trace of $A^n$ can be larger than
the number $F_n$ of period $n$ points of $h$.  A combinatorial argument of Manning  \cite{Man71} can be used to correct the count.  (For a homological approach, see   \cite{Fri87}.) One constructs additional shifts of finite type $\sigma_i,\ i=1, \ldots, m$, with adjacency matrices $A_i$ and signs
$\epsilon_i \in \{ -1, 1\}$ such that for each $n$,
\begin{equation}\label{eq4.2}F_n = \sum_{i=0}^m \epsilon_i \cdot tr A_i^n.\end{equation}
Let $p(t)=\prod(t-\lambda_j)$ where $\lambda_j$ ranges over the distinct non-zero roots  of the characteristic polynomials of $A_0, \ldots, A_m$. Then $F_n$ is linearly recurrent with minimal polynomial
$p(t)$. 
\end{proof}

The invariant foliations are {\it orientable} if all of the intersections of any transverse arc with any leaf are in the same direction. For surfaces with connected boundary, the condition is equivalent to the requirement that the degree of every interior prong singularity is even.

The pseudo-Anosov monodromy of the $(-2,3,7)$-pretzel knot has orientable
invariant foliations (see remark following Theorem 9.7 of \cite{Lei03}). 
For it, $\max_k \GR^{(k)}({\bf F})= \GR^{(1)}({\bf F})$ is equal to the Mahler measure of Lehmer's polynomial. More generally, we have

\begin{theorem}\label{thm4.4} Assume that $h: S \to S$ is a pseudo-Anosov homeomorphism with oriented invariant foliations. Then for any positive integer $k$, 
$$\GR^{(k)}({\bf F}) = \GR^{(k)}({\bf L}).$$\end{theorem}

Combined with Proposition \ref{prop4.1} this gives: 

\begin{corollary}\label{cor4.5}  Assume that $h: S \to S$ is a pseudo-Anosov homeomorphism with oriented invariant foliations and $f(t) = \det (h_* - t I)$. Then 
$$\max_k \GR^{(k)}({\bf F}) \le M(f),$$
with equality if $f$ has no repeated roots.\end{corollary}

\begin{proof}[Proof of Theorem~\textup{4$\cdot$4.}] For any $n\ge 1$, all non-prong fixed points of $h^n$ have the same index, $+1$ or $-1$. One way to see this is to choose such a  fixed point $p$, and orient the leaf $L^u_p$ of  ${\cal F}^u$ containing $p$. The index at $p$ is $-1$ if $h$ preserves the orientation of $L^u_p$; otherwise it is $+1$ (see Proposition 5.7 of \cite{Fra82}). Let $q$ be any other non-prong fixed point of $h^n$, and let $L^u_q$ be the leaf of  ${\cal F}^u$ containing it. Since $L^u_p$ is dense in $S$, it passes arbitrarily close to $q$, always in the same direction, since ${\cal F}^u$ is oriented. By continuity, the orientations of $L_q^u$ and $L_q^u$ are both preserved or else both reversed. Thus the index of $q$ agrees with that of $p$.

If $h$ reverses orientations of unstable leaves, then replace each $F_n$ by $(-1)^{n-1} F_n$.  Otherwise, replace each $F_n$ by
$-F_n$.  In either case, the generalized growth rates are unchanged.  Any remaining difference between $F_n$  and $L_n$ is due to prong singularities, and hence is bounded.  Proposition \ref{prop2.3} completes the proof. \end{proof}

\begin{remark}\label{rmk4.6} Theorem 3.3 of  \cite{Ry99} states that for any pseudo-Anosov homeomorphism $h$ with orientable foliations, the eigenvalues of $h_*$ are the same as those of the matrix (\ref{eq4.1}) including multiplicity, with the possible exception of some zeros and roots of unity. Corollary \ref{cor4.5} above asserts that in the case that $h$ has oriented invariant foliations, the generalized growth rates of the sequence $\{{\rm trace}\ A^n\}$ are equal to those of ${\bf F}$. \end{remark}

\section{Free group endomorphisms}\label{sec5}  Let $G$ be any group with a finite set of generators $g_1, \ldots, g_m$. The {\it length} $|g|$ of any $1\ne g\in G$ is the minimum length of a word in $g_1^{\pm 1}, \ldots, g_m^{\pm 1}$ expressing $g$. If $\phi$ is any endomorphism of  $G$, then the {\it growth rate} $\GR(\phi)$ is 
$$\GR(\phi) = \max_{1\le i \le m}\limsup_{n\to \infty} |\phi^n (g_i)|^{1/n}.$$
It is  independent of the generator set  \cite{Bow78}. 

When $f:S\to S$ is a pseudo-Anosov homeomorphism, the growth rate $\GR(f_\sharp)$ of the induced automorphism $f_\sharp$ of $G = \pi_1(S)$ contains a wealth of information about the dynamics of $f$. For example, $\GR(f_\sharp)$ coincides with  the first growth rate  $\GR^{(1)}({\bf F})$, and it is an upper bound for the modulus of every eigenvalue of $f_*: H_1(S; {\mathbb R}) \to H_1(S; {\mathbb R})$ (see  \cite{FS79} and  \cite{Bow78}). Hence
\begin{equation}\label{eq5.1}\GR(f_\sharp) = \GR^{(1)}({\bf F}) \ge \GR^{(1)}({\bf L}).\end{equation}

Topological entropy $h_{\rm top}(f)$ is a measure of complexity defined for any continuous map $f: X \to X$ of a compact metric space. The entropy of any periodic map is zero, while that of a pseudo-Anosov map is strictly positive. Although topological entropy is difficult to compute in general, when $X$ is a surface and $f$ is pseudo-Anosov, $h_{\rm top}(f)$ can be computed as $\log \GR(f_\sharp)$  \cite{FS79}.

\begin{definition} Let $G$ be a group generated by $g_1, \ldots, g_m$, and let $\phi: G \to G$ be an endomorphism. The {\it $k$th generalized growth rate}
of $\phi$  (with respect to the given generating set) is 
$$\GR^{(k)}(\phi) = \max_{1\le i \le m}\GR^{(k)}( |\phi^n (g_i)|).$$\end{definition}

\begin{example}\label{ex5.2} In general these growth rates may depend on the generating set when $k>1$. Consider the free group $F$ of rank two generated by $x, y$ and the endomorphism $\phi: F \to F$ sending
$x$ to $x^3$ and $y$ to $y^2$. It  is easy to see that 
$|\phi^n(x)| = 3^n$,\  $ |\phi^n(y)| = 2^n$ and hence $\GR^{(1)}(\phi) = 3$ while $\GR^{(k)}(\phi) = 0$ for $k>1$. 

On the other hand, if we replace $x, y$ with generators $z, y$, where $z = xy$, then 
$\phi^n(z)=  x^{3^n} y^{2^n} = (zy^{-1})^{3^n -1} z y^{2^n -1}$ and 
$|\phi^n(z)| = 2\cdot 3^n + 2^n -2$. The  latter is linearly recurrent with minimal polynomial $(t-1)(t-2)(t-3)$ and hence $\GR^{(1)}(\phi) =3, 
\GR^{(2)}(\phi) = \GR^{(3)}(\phi) = 6$ while $ \GR^{(k)}(\phi) =0$ for $k >3$.
\end{example}

\begin{remark}\label{altdef} If we define $\GR^{(k)}(\phi)$ instead to be $\GR^{(k)}( \sum_{i=1}^m|\phi^n (g_i)|)$ then the two generating sets in the above example give the same maximum growth rate.  We do not know if this holds in general.  It is easily seen that this definition reduces to the usual growth rate $\GR(\phi)$ when $k=1$.
\end{remark}

An important case of Lehmer's question can be viewed as a question about free group endomorphisms. 

A square matrix $A$ with nonnegative real entries is {\it primitive} if there is an $N>0$ such that all the entries of $A^N$ are positive. The characteristic polynomial of any primitive integer matrix is a monic integral polynomial such that one root is a positive real number that is strictly larger than the modulus of the remaining roots. Such a polynomial is said to be a {\it Perron polynomial}, and the dominant root is called a {\it Perron number}. For background about Perron--Frobenius theory, the reader might consult  \cite{Mi88}. 

It follows from equation (\ref{eq4.2}) that the minimal polynomial $p(t)$ of ${\bf F}$ is a product of Perron polynomials. 

\begin{theorem}\label{perron} Given any  irreducible Perron polynomial $p(t)$ there exists an endomorphism $\phi$ of a finitely generated free group such that for any generator $g_i$, the sequence of word lengths  $|\phi^n (g_i)|$ is linearly recurrent with minimal polynomial equal to $p(t)$ times a product of cyclotomic polynomials. Thus the maximum growth rate of the sequence is $M(p)$.\end{theorem}

\begin{proof}  For $\Lambda=(\lambda_1, \lambda_2, \ldots, \lambda_d)$  a $d$-tuple of complex numbers, we define ${\rm tr}(\Lambda^n)=\sum_{i=1}^d \lambda_i^n$ and ${\rm tr}_n(\Lambda)=\sum_{k|n}\mu(n/k){\rm tr}(\Lambda^k)$ where $\mu$ is the Moebius function.  (The second quantity is called the $n$th {\it net trace} of $\Lambda$.)
The main theorem of  \cite{KOR00} states that  there exist a primitive integer matrix $A$ and a non-negative integer $\ell$ such that ${\rm char}(A)=t^\ell\prod_{i=1}^d( t-\lambda_i)$ if and only if: 
\begin{enumerate}
\item the polynomial $\prod_{i=1}^d (t-\lambda_i )$ has integer coefficients;
\item some $\lambda_i$ is real and $\lambda_i > |\lambda_j|$ for $i\neq j$;
\item ${\rm tr}_n(\Lambda) \ge 0$ for all $n\ge 1$. 
\end{enumerate}
Let  $p(t)=\prod_{i=1}^d( t-\lambda_i )$ be a Perron polynomial. Since one root dominates the moduli of the other roots, the net trace condition is satisfied for sufficiently large $n$. We observe that the $k$-tuple of $k$th roots of unity has $n$th net trace equal to $k$ for $n=k$ and 0  otherwise.  Thus we can multiply $p(t)$ by a suitable product $\Phi(t)$ of cyclotomic polynomials to obtain a polynomial $p(t)\Phi(t)$ for which the tuple of roots satisfies  the net trace condition for all $n\ge 1$. 
Then using  \cite{KOR00},  we can find a primitive integer $m\times m$ matrix  
$A = (a_{i,j})$ with characteristic polynomial equal to $t^\ell p(t)\Phi(t)$ for some $\ell$. 

Let $F_m$ be the free group generated by $g_1, \ldots, g_m$. Define $\phi: F_m\to F_m$ to be the endomorphism mapping $g_i$ to $g_1^{a_{i,1}}g_2^{a_{i,2}}\cdots g_m^{a_{i,m}}$. 

Only nonnegative exponents occur in any reduced word representing $\phi^n(g_i)$, for any $n\ge 1$ and  $i=1, \ldots, m$. Hence the number of occurrences of $g_j$ in $\phi^n(g_i)$ is equal to $(A^n)_{i,j}$. The Cayley--Hamilton theorem implies that for fixed $i$ and  $j$, the sequence $(A^n)_{i,j}$ satisfies a
linear recurrence with characteristic polynomial ${\rm char}(A)$. Hence for any $i$, so does the sequence of word lengths  $|\phi^n (g_i)|$, which is equal to $\sum_j (A^n)_{i,j}$. Since $A$ is primitive, the sequence grows exponentially at a rate given by the dominant eigenvalue.  Hence the minimal polynomial of the sequence is equal to the Perron polynomial $p(t)$ times some factor of $\Phi(t)$.  \end{proof}




If $\ell =0$, then $\phi$ is injective since its image abelianizes to a free abelian group of rank $m$. In such a case, we cannot expect $\phi: F \to F$ in Theorem 5.3 to be induced by an automorphism of an orientable surface, since in such case the set of roots of $p(t)$ must be closed under inversion (see  \cite{St82}). Instead we ask:

\begin{question} Can $\phi: F \to F$ in Theorem \ref{perron} be chosen to be a free group {\sl automorphism}? \end{question}

It is well known that given any matrix $A \in GL_m(\Z)$, there exists a free group automorphism $\phi \in {\rm Aut}(F_m)$ that abelianizes to $A$ (see  \cite{LS77}). We say an endomorphism of the free group is {\sl positive} if there is a generating set for which the image of every generator is a positive word, that is, a product of positive powers of generators.  

\begin{question}\label{q5.5} Suppose that $A \in {\rm GL}_m(\Z)$ has nonnegative entries. Does there exist a positive automorphism of the free group $F_m$ that abelianizes to $A$? \end{question}

We conclude the section with a partial answer to Question \ref{q5.5}.

\begin{proposition}\label{prop5.6}  Suppose that $A \in {\rm GL}_2(\Z)$ has nonnegative entries. Then there is a positive automorphism of the free group $F$ of rank 2 that abelianizes to $A$.\end{proposition}

\begin{proof} Write 
$$A=\begin{pmatrix}p_0 \ & p_1\\ q_0 & q_1\end{pmatrix}.$$
We may assume (exchanging the two generators if necessary) that either $p_0-p_1$ or $q_0-q_1$ is positive; the other is necessarily nonnegative. 

Inductively, for $k\ge 0$ we write
$$\begin{pmatrix}p_k\\ q_k\end{pmatrix}= d_{k+1} \begin{pmatrix}p_{k+1}\\ q_{k+1}\end{pmatrix}+ \begin{pmatrix}p_{k+2}\\ q_{k+2}\end{pmatrix},$$
where $d_{k+1}$ is chosen to be the largest integer with $p_{k+2},
q_{k+2}$ both nonnegative. We stop when $n=k+2$ satisfies 
$p_n=0$ or $q_n=0$. We see easily that each matrix
$$\begin{pmatrix}p_k \ & p_{k+1}\\ q_k & q_{k+1}\end{pmatrix}$$ 
is in ${\rm GL}_2(\Z)$
with $p_k-p_{k+1}$ or $q_k-q_{k+1}$ positive.  In particular, if $p_n=0$ then $q_n=p_{n-1}=1$, while if $q_n=0$ then $p_n=q_{n-1}=1$.

Let $F_2$ be the free group of rank 2 generated by $a,b$. Set
$u_0 = a^{p_n}b^{q_n},\ u_1= a^{p_{n-1}} b^{q_{n-1}}$ (these generate $F_2$), and 
$u_{k+1} = u_k^{d_{n-k}}u_{k-1}$ for $k=1, \ldots, n-1$. Define
$\phi: F_2 \to F_2$ by $\phi(a) = u_n,\ \phi(b) = u_{n-1}$. We see by descent that $u_k$ is in the image of $\phi$ for all $k$, and so $\phi$ is onto. Any surjective endomorphism of $F_2$ is an automorphism.\end{proof}


\section{Braids and Lehmer's question}\label{sec6}  An $n$-braid can be regarded as an orientation-preserving homeomorphism of the punctured disk $D_n= D^2\setminus \{ n\ {\rm points}\}$, two $n$-braids regarded as the same if they are isotopic rel boundary. Alternatively, one can regard an $n$-braid as an isotopy class of a map of the pair $(D^2, \{n\ {\rm points}\})$. The latter point of view enables us to recover the usual geometric picture of an $n$-braid by tracing the image of the set of $n$ points as an isotopy rel boundary is performed from the map of $D^2$ to the identity. 

The collection of all $n$-braids forms a group $B_n$ under concatenation. As usual, we let $\s_i$ denote the final map of an isotopy of $D_n$ that exchanges the $i$th and $i+1$st puncture, dragging the $i$th puncture around the $i+1$st in the clockwise direction. Then $B_n$ is generated by $\s_1, \ldots, \s_{n-1}$. The center of $B_n$ is cyclic, generated by a full twist $\D_n^2 = (\s_{n-1}\s_{n-2}\cdots \s_1)^n$. 

The Thurston--Nielsen classification of surface homeomorphisms applies to $n$-braids. One calls $\beta\in B_n$ periodic, reducible or pseudo-Anosov if it can be represented by a homeomorphism of $D_n$ with that property. 

Given an $n$-braid $\beta$, its {\it braid entropy} $h_\beta$ is the infimum
${\rm inf}\{h_{\rm top}(f)\mid f\ {\rm represents}\ \beta \}$. Braid entropy is a conjugacy invariant. Moreover, the braid entropy of $\beta\in B_n$ is equal to that of 
$\D_n^{2k}\beta$, for any integer $k$. If $f:D_n \to D_n$ is any representative of $\beta$, then $h_\beta$ can be computed as $\log \GR(f_\sharp)$. 

The {\it reduced Burau representation} $\beta \mapsto \tilde B_\beta$ maps $B_n$ to ${\rm GL}(n-1, \Z[t^{\pm 1}])$ via 
$$ \tilde B_{\s_1}=\left(
\begin{array}{c@{\enskip}c@{\enskip}c}
-t&0&\\ 1 &1&\\&&I
\end{array}\right),\quad \tilde B_{\s_i}=\left(
\begin{array}{c@{\enskip}c@{\enskip}c@{\enskip}c@{\enskip}c}
I&&&&\\ &1&t&0&\\ &0&-t&0&\\ &0&1&1&\\&&&&I
\end{array}\right),\quad \tilde B_{\s_{n-1}}=\left(
\begin{array}{c@{\enskip}c@{\enskip}c}
I&&\\&-t&0\\&1 &1
\end{array}\right),$$
where $I$ denotes an identity matrix of the appropriate size. 

The closure $\hat \beta$ of an $n$-braid is the oriented link obtained  by joining the top and bottom of each strand of $\beta$ without introducing additional crossings. 
The strands of $\beta$ are oriented coherently, from top to bottom. The following 
is well known (see  \cite{Bi74}). 

\begin{proposition}\label{prop6.1} The matrix $\tilde B_\beta -I$ is a Jacobian of the link $\hat \beta$. Moreover, ${\rm det}(\tilde B_\beta-I)$ is equal to $(1+t+\cdots +t^{n-1}) \D_{\rm red}(\hat \beta)(t)$, where $\D_{\rm red}$ is the reduced Alexander polynomial. \end{proposition}

As a consequence of results of section \ref{sec3} we have: 

\begin{proposition}\label{prop6.2} Lehmer's question is equivalent to the following.
For any $\epsilon >0$, does there exist a hyperbolic $n$-braid $\beta$ such that $1<M({\rm det}(B_\beta-I))<1+\epsilon$?\end{proposition}

The example below suggests that the relationship between Lehmer's question and surface automorphisms is significant. We see that for $n=3,4$ and $5$, there exists a minimal entropy pseudo-Anosov $n$-braid $f: D_n \to D_n$ that closes to the 
$(-2,3,7)$-pretzel knot. 

In each case, the mapping torus of $f$ is a hyperbolic $3$-manifold $M$ with $H_1M \cong \Z^2$. There are many possible
homomorphisms $\pi_1M \to \Z$ with finitely generated kernel, homomorphisms that one might call fibering directions. (Fibering directions can be viewed as points of fibered faces of a polyhedral unit ball of the Thurston norm. See  \cite{Mc03}.)
One such direction yields the punctured disk $D_n$ as fiber and $f$ as monodromy. Other fibering directions result in pseudo-Anosov monodromies of possibly higher genus surfaces. In each of the three cases below, there exists a fibering direction for which the corresponding monodromy, capped off along the boundary of $D$, is that of the $(-2,3,7)$-pretzel knot. 
 
\begin{example}\label{ex6.3}
(1) The braid $\s_1 \s_2^{-1}$ is minimal among all pseudo-Anosov $3$-braids in the forcing partial order  \cite{Mat85, Han97}. Hence it attains the minimum braid entropy (=$\log 2.61803...$) among all pseudo-Anosov $3$-braids. We note here that the closure of $\s_1 \s_2^{-1}\D_3^4$, which has the same braid entropy, is the $(-2,3,7)$-pretzel knot. 

(2) The braid $\s_3\s_2\s_1^{-1}$ attains the minimum braid entropy (=$\log 2.29663...$) among all pseudo-Anosov $4$-braids  \cite{HS05}. The closure of $\s_3\s_2\s_1^{-1} \D_4^2$ is the $(-2,3,7)$-pretzel knot. 

(3) The braid $\s_1\s_2\s_3\s_4\s_1\s_2$ attains the minimum braid entropy (=$\log 1.72208...$) among all pseudo-Anosov $5$-braids  \cite{HS05}.  The  closure
of $\s_1\s_2\s_3\s_4\s_1\s_2 \D_5^2$ is again the $(-2,3,7)$-pretzel knot. 
\end{example}

We expect the mapping torus of a pseudo-Anosov homeomorphims with small entropy to be a hyperbolic $3$-manifold that is small in the sense of  \cite{CDW98}; that is, composed of relatively few ideal tetrahedra. The examples above do not disappoint. The first and third have 4 tetrahedra, while the second has 5. Their volumes are 4.0597..., 3.6638... and 
4.8511..., respectively.  Calculations were done with SnapPea. We note that the volume of the third example is the same as that of the Whitehead link complement, and is the smallest known volume of any hyperbolic 2-component link.


\end{document}